# DISCUSSION OF "LEAST ANGLE REGRESSION" BY EFRON ET AL.

By Robert A. Stine

*University of Pennsylvania*

I have enjoyed reading the work of each of these authors over the years, so it is a real pleasure to have this opportunity to contribute to the discussion of this collaboration. The geometry of LARS furnishes an elegant bridge between the Lasso and Stagewise regression, methods that I would not have suspected to be so related. Toward my own interests, LARS offers a rather different way to construct a regression model by gradually blending predictors rather than using a predictor all at once. I feel that the problem of "automatic feature generation" (proposing predictors to consider in a model) is a current challenge in building regression models that can compete with those from computer science, and LARS suggests a new approach to this task. In the examples of Efron, Hastie, Johnstone and Tibshirani (EHJT) (particularly that summarized in their Figure 5), LARS produces models with smaller predictive error than the old workhorse, stepwise regression. Furthermore, as an added bonus, the code supplied by the authors runs faster for me than the `step` routine for stepwise regression supplied with $R$, the generic version of S-PLUS that I use.

My discussion focuses on the use of $C_p$ to choose the number of predictors. The bootstrap simulations in EHJT show that LARS reaches higher levels of "proportion explained" than stepwise regression. Furthermore, the goodness-of-fit obtained by LARS remains high over a wide range of models, in sharp contrast to the narrow peak produced by stepwise selection. Because the cost of overfitting with LARS appears less severe than with stepwise, LARS would seem to have a clear advantage in this respect. Even if we do overfit, the fit of LARS degrades only slightly. The issue becomes learning how much LARS overfits, particularly in situations with many potential predictors ($m$ as large as or larger than $n$).

To investigate the model-selection aspects of LARS further, I compared LARS to stepwise regression using a "reversed" five-fold cross-validation. The cross-validation is reversed in the sense that I estimate the models on one fold (88 observations) and then predict the other four. This division with more set aside for validation than used in estimation offers a better comparison of models. For example, Shao (1993) shows that one needs to let the proportion used for validation grow large in order to get cross validation to find the right model. This reversed design also adds a further benefit of making the variable selection harder. The quadratic model fitted to the diabetes data in EHJT selects from $m = 64$ predictors using a sample of $n = 442$ cases, or about 7 cases per predictor. Reversed cross-validation is closer to a typical data-mining problem. With only one fold of 88 observations to train the model, observation noise obscures subtle predictors. Also, only a few degrees of freedom remain to estimate the error variance $\overline{\sigma}^2$ that appears in $C_p$ [equation (4.5)]. Because I also wanted to see what happens when $m > n$, I repeated the cross-validation with 5 additional possible predictors added to the 10 in the diabetes data. These 5 spurious predictors were simulated i.i.d. Gaussian noise; one can think of them as extraneous predictors that one might encounter when working with an energetic, creative colleague who suggests many ideas to explore. With these 15 base predictors, the search must consider $m = 15 + \binom{15}{2} + 14 = 134$ possible predictors.

Here are a few details of the cross-validation. To obtain the stepwise results, I ran forward stepwise using the hard threshold $2 \log m$, which is also known as the risk inflation criterion (RIC) [Donoho and Johnstone (1994) and Foster and George (1994)]. One begins with the most significant predictor. If the squared $t$-statistic for this predictor, say $t^2_{(1)}$, is less than the threshold $2 \log m$, then the search stops, leaving us with the "null" model that consists of just an intercept. If instead $t^2_{(1)} \geq 2 \log m$, the associated predictor, say $X_{(1)}$, joins the model and the search moves on to the next predictor. The second predictor $X_{(2)}$ joins the model if $t^2_{(2)} \geq 2 \log m$; otherwise the search stops with the one-predictor model. The search continues until reaching a predictor whose $t$-statistic fails this test, $t^2_{(q+1)} < 2 \log m$, leaving a model with $q$ predictors. To obtain the results for LARS, I picked the order of the fit by minimizing $C_p$. Unlike the forward, sequential stepwise search, LARS globally searches a collection of models up to some large size, seeking the model which has the smallest $C_p$. I set the maximum model size to 50 (for $m = 64$) or 64 (for $m = 134$). In either case, the model is chosen from the collection of linear and quadratic effects in the 10 or 15 basic predictors. Neither search enforces the principle of marginality; an interaction can join the model without adding main effects.

I repeated the five-fold cross validation 20 times, each time randomly partitioning the 442 cases into 5 folds. This produces 100 stepwise and LARS fits. For each of these, I computed the square root of the out-of-sample mean square error (MSE) when the model fit on one fold was used to predict the held-back $354 \ [= 4(88) + 2]$ observations. I also saved the size $q$ for each fit.

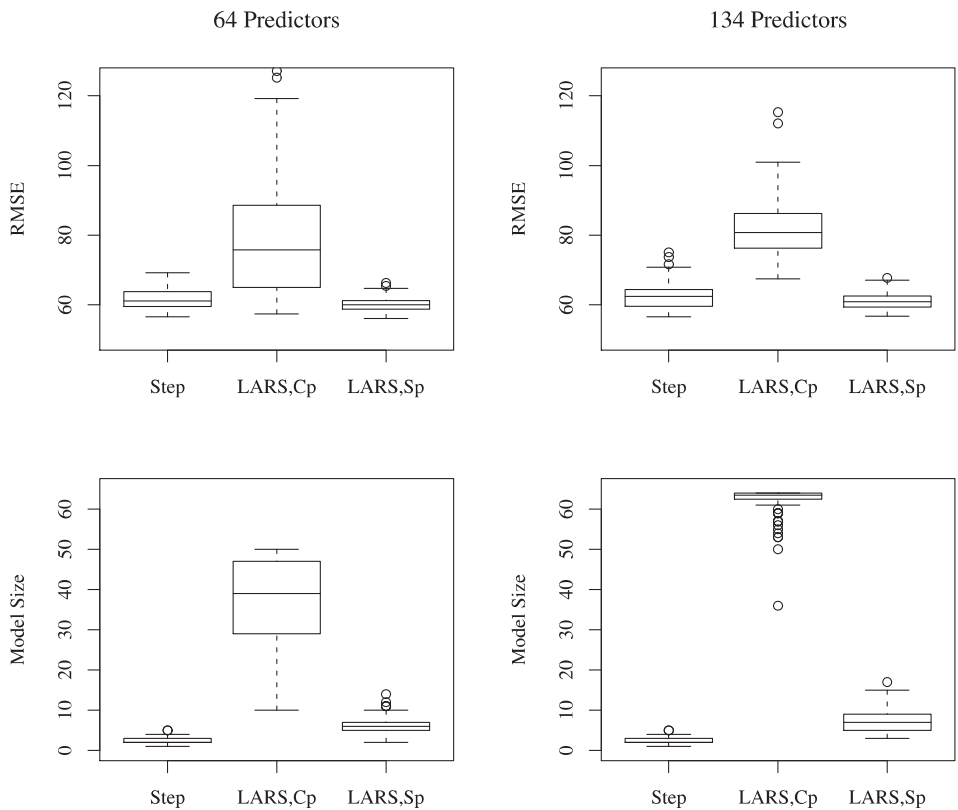

FIG. 1. *Five-fold cross-validation of the prediction error and size of stepwise regression and LARS when fitting models to a collection of* 64 (left) *or* 134 *predictors* (right). *LARS fits chosen by $C_p$ overfit and have larger RMSE than stepwise; with $C_p$ replaced by the alternative criterion $S_p$ defined in* (3), *the LARS fits become more parsimonious with smaller RMSE. The random splitting into estimation and validation samples was repeated* 20 *times, for a total of* 100 *stepwise and LARS fits.*

Figure 1 summarizes the results of the cross-validation. The comparison boxplots on the left compare the square root of the MSE (top) and selected model order (bottom) of stepwise to LARS when picking from $m = 64$ candidate predictors; those on the right summarize what happens with $m = 134$. When choosing from among 64 predictors, the median size of a LARS model identified by $C_p$ is 39. The median stepwise model has but 2 predictors. (I will describe the $S_p$ criterion further below.) The effects of overfitting on the prediction error of LARS are clear: LARS has higher RMSE than stepwise. The median RMSE for stepwise is near 62. For LARS, the median RMSE is larger, about 78. Although the predictive accuracy of LARS declines more slowly than that of stepwise when it overfits (imagine the fit of stepwise with 39 predictors), LARS overfits by enough in this case that it predicts worse than the far more parsimonious stepwise models. With more predictors ($m = 134$), the boxplots on the right of Figure 1 show that $C_p$ tends to pick the largest possible model—here a model with 64 predictors.

Why does LARS overfit? As usual with variable selection in regression, it is simpler to try to understand when thinking about the utopian orthogonal regression with known $\sigma^2$. Assume, as in Lemma 1 of EHJT, that the predictors $X_j$ are the columns of an identity matrix, $X_j = e_j = (0, \ldots, 0, 1_j, 0, \ldots, 0)$. Assume also that we know $\sigma^2 = 1$ and use it in place of the troublesome $\overline{\sigma}^2$ in $C_p$, so that for this discussion

$$C_p = \text{RSS}(p) - n + 2p. \tag{1}$$

To define RSS($p$) in this context, denote the ordered values of the response as

$$Y_{(1)}^2 > Y_{(2)}^2 > \cdots > Y_{(n)}^2.$$

The soft thresholding summarized in Lemma 1 of EHJT implies that the residual sum-of-squares of LARS with $q$ predictors is

$$\text{RSS}(q) = (q+1)Y_{(q+1)}^2 + \sum_{j=q+2}^{n} Y_{(j)}^2.$$

Consequently, the drop in $C_p$ when going from the model with $q$ to the model with $q+1$ predictors is

$$C_q - C_{q+1} = (q+1)\,d_q - 2,$$

with

$$d_q = Y_{(q+1)}^2 - Y_{(q+2)}^2;$$

$C_p$ adds $X_{q+1}$ to the model if $C_q - C_{q+1} > 0$.

This use of $C_p$ works well for the orthogonal "null" model, but overfits when the model includes much signal. Figure 2 shows a graph of the mean and standard deviation of $\text{RSS}(q) - \text{RSS}(0) + 2q$ for an orthogonal model with $n = m = 100$ and $Y_i \stackrel{\text{i.i.d.}}{\sim} N(0,1)$. I subtracted RSS(0) rather than $n$ to reduce the variation in the simulation. Figure 3 gives a rather different impression. The simulation is identical except that the data have some signal. Now, $EY_i = 3$ for $i = 1, \ldots, 5$. The remaining 95 observations are $N(0,1)$. The "true" model has only 5 nonzero components, but the minimal expected $C_p$ falls near 20.

This stylized example suggests an explanation for the overfitting—as well as motivates a way to avoid some of it. Consider the change in RSS for a null model when adding the sixth predictor. For this step, $\text{RSS}(5) - \text{RSS}(6) = 6(Y_{(6)}^2 - Y_{(7)}^2)$. Even though we multiply the difference between the squares by 6, adjacent order statistics become closer when removed from the extremes, and $C_p$ tends to increase as shown in Figure 2. The situation changes when signal is present. First, the five observations with mean 3 are likely to be the first five ordered observations. So, their spacing is likely to be larger because their order is determined by a sample of five normals; $C_p$ adds these. When reaching the noise, the difference $Y_{(6)}^2 - Y_{(7)}^2$ now behaves like the difference between the *first two* squared order statistics in an i.i.d. sample of 95 standard normals. Consequently, this comparison involves the gap between the most extreme order statistics rather than those from within the sample, and as a result $C_p$ drops to indicate a larger model.

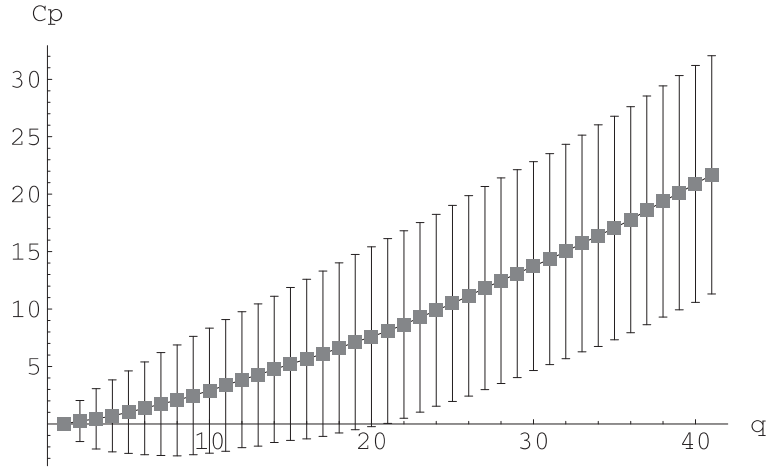

Fig. 2. *A simulation of $C_p$ for LARS applied to orthogonal, normal data with no signal correctly identifies the null model. These results are from a simulation with 1000 replications, each consisting of a sample of 100 i.i.d. standard normal observations. Error bars indicate $\pm 1$ standard deviation.*

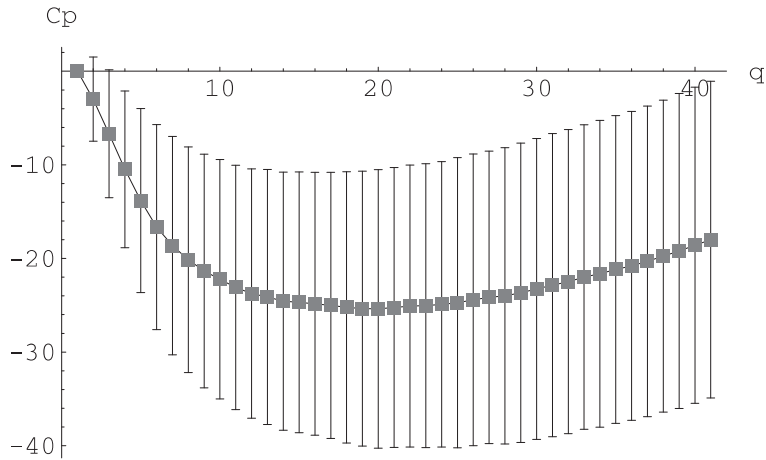

Fig. 3. *A simulation of $C_p$ for LARS applied to orthogonal, normal data with signal present overfits. Results are from a simulation with 1000 replications, each consisting of 5 observations with mean 3 combined with a sample of 95 i.i.d. standard normal observations. Error bars indicate $\pm 1$ standard deviation.*

This explanation of the overfitting suggests a simple alternative to $C_p$ that leads to smaller LARS models. The idea is to compare the decreasing residual sum of squares $\text{RSS}(q)$ to what is expected under a model that has fitted some signal *and* some noise. Since overfitting seems to have relatively benign effects on LARS, one does not want to take the hard-thresholding approach; my colleague Dean Foster suggested that the criterion might do better by assuming that some of the predictors already in the model are really noise. The criterion $S_p$ suggested here adopts this notion. The form of $S_p$ relies on approximations for normal order statistics commonly used in variable selection, particularly adaptive methods [Benjamini and Hochberg (1995) and Foster and Stine (1996)]. These approximate the size of the $j$th normal order statistic in a sample of $n$ with $\sqrt{2\log(n/j)}$. To motivate the form of the $S_p$ criterion, I return to the orthogonal situation and consider what happens when deciding whether to increase the size of the model from $q$ to $q+1$ predictors. If I know that $k$ of the already included $q$ predictors represent signal and the rest of the predictors are noise, then $d_q = Y_{(q+1)}^2 - Y_{(q+2)}^2$ is about

$$(2) \qquad 2\log\frac{m-k}{q+1-k} - 2\log\frac{m-k}{q+2-k}.$$

Since I do not know $k$, I will just set $k = q/2$ (i.e., assume that half of those already in the model are noise) and approximate $d_q$ as

$$\delta(q) = 2\log\frac{q/2+2}{q/2+1}.$$

[Define $\delta(0) = 2\log 2$ and $\delta(1) = 2\log 3/2$.] This approximation suggests choosing the model that minimizes

$$(3) \qquad S_q = \text{RSS}(q) + \hat{\sigma}^2 \sum_{j=1}^{q} j\delta(j),$$

where $\hat{\sigma}^2$ represents an "honest" estimate of $\sigma^2$ that avoids selection bias. The $S_p$ criterion, like $C_p$, penalizes the residual sum-of-squares, but uses a different penalty.

The results for LARS with this criterion define the third set of boxplots in Figure 1. To avoid selection bias in the estimate of $\sigma^2$, I used a two-step procedure. First, fit a forward stepwise regression using hard thresholding. Second, use the estimated error variance from this stepwise fit as $\hat{\sigma}^2$ in $S_p$ and proceed with LARS. Because hard thresholding avoids overfitting in the stepwise regression, the resulting estimator $\hat{\sigma}^2$ is probably conservative—but this is just what is needed when modeling data with an excess of possible predictors. If the variance estimate from the largest LARS model is used instead, the $S_p$ criterion also overfits (though not so much as $C_p$). Returning to Figure 1, the combination of LARS with $S_p$ obtains the smallest typical MSE with both $m = 64$ and 134 predictors. In either case, LARS includes more predictors than the parsimonious stepwise fits obtained by hard thresholding.

These results lead to more questions. What are the risk properties of the LARS predictor chosen by $C_p$ or $S_p$? How is it that the number of possible predictors $m$ does not appear in either criterion? This definition of $S_p$ simply supposes half of the included predictors are noise; why half? What is a better way to set $k$ in (2)? Finally, that the combination of LARS with either $C_p$ or $S_p$ has less MSE than stepwise when predicting diabetes is hardly convincing that such a pairing would do well in other applications. Statistics would be well served by having a repository of test problems comparable to those held at UC Irvine for judging machine learning algorithms [Blake and Merz (1998)].

DEPARTMENT OF STATISTICS
THE WHARTON SCHOOL
UNIVERSITY OF PENNSYLVANIA
PHILADELPHIA, PENNSYLVANIA 19104-6340
USA
E-MAIL: stine@wharton.upenn.edu